

\def\TheMagstep{\magstep1}	
 \def\PaperSize{AFour}		
\let\:=\colon \def\IP{{\bf P}} \let\x=\times 
\def\O{{\cal O}} \def\I{{\cal I}}  
 \def\N{{\cal N}}
\def\and{\hbox{ and }} \def\for{\hbox{ for }}
\let\?=\overline
\def\Hilb{\mathop{\rm Hilb}\nolimits}

\parskip=0pt plus 1.75pt
\parindent10pt
\hsize26pc
\vsize=44\normalbaselineskip 

\abovedisplayskip6pt plus6pt minus0.25pt
\belowdisplayskip6pt plus6pt minus0.25pt

\def\TRUE{TRUE}
\ifx\DoublepageOutput\TRUE \def\TheMagstep{\magstep0} \fi
\mag=\TheMagstep

\newskip\vadjustskip \vadjustskip=0.5\normalbaselineskip
\def\centertext
 {\hoffset=\pgwidth \advance\hoffset-\hsize
  \advance\hoffset-2truein \divide\hoffset by 2\relax
  \voffset=\pgheight \advance\voffset-\vsize
  \advance\voffset-2truein \divide\voffset by 2\relax
  \advance\voffset\vadjustskip
 }
\newdimen\pgwidth\newdimen\pgheight
\def\letter{letter}\def\AFour{AFour}
\ifx\PaperSize\letter
 \pgwidth=8.5truein \pgheight=11truein
 \message{- Got a paper size of letter.  }\centertext
\fi
\ifx\PaperSize\AFour
 \pgwidth=210truemm \pgheight=297truemm
 \message{- Got a paper size of AFour.  }\centertext
\fi

 \newdimen\fullhsize \newbox\leftcolumn
 \def\fulline{\hbox to \fullhsize}
\def\doublepageoutput
{\let\lr=L
 \output={\if L\lr
          \global\setbox\leftcolumn=\columnbox \global\let\lr=R%
        \else \doubleformat \global\let\lr=L\fi
        \ifnum\outputpenalty>-20000 \else\dosupereject\fi}%
 \def\doubleformat{\shipout\vbox{%
        \fulline{\hfil\hfil\box\leftcolumn\hfil\columnbox\hfil\hfil}%
}%
 }%
 \def\columnbox{\vbox
   {\makeheadline\pagebody\makefootline\advancepageno}%
   }
 \fullhsize=\pgheight \hoffset=-1truein
 \voffset=\pgwidth \advance\voffset-\vsize
  \advance\voffset-2truein \divide\voffset by 2
  \advance\voffset\vadjustskip
\let\firstheadline=\hfil

\null\vfill\nopagenumbers\eject\pageno=1\relax 
}

\catcode`\@=11
\def\vfootnote#1{\insert\footins\bgroup
 \eightpoint 
 \interlinepenalty\interfootnotelinepenalty
  \splittopskip\ht\strutbox 
  \splitmaxdepth\dp\strutbox \floatingpenalty\@MM
  \leftskip\z@skip \rightskip\z@skip \spaceskip\z@skip \xspaceskip\z@skip
  \textindent{#1}\footstrut\futurelet\next\fo@t}

\newif\ifdates	\newif\ifdateonpageone
\def\today{\ifcase\month\or
 January\or February\or March\or April\or May\or June\or
 July\or August\or September\or October\or November\or December\fi
 \space\number\day, \number\year}

\nopagenumbers
\headline={%
 \eightpoint
  \ifnum\pageno=1\firstheadline
  \else
    \ifodd\pageno\oddheadline
    \else\evenheadline\fi
  \fi
}

 \def\firstheadline{\ifdateonpageone\rm\today\fi \hfill}
 \def\oddheadline{\rm
  \ifdates\rlap{March 20, 1996}\fi
  \hfil\shorttitle\hfil\llap{\the\pageno}}
 \def\evenheadline{\rm\rlap{\number\pageno}\hfil
 	\author\hfil
	\ifdates\llap{alg-geom/9601024}\fi}

 \newcount\sectno \sectno=0
 \newskip\sectskipamount \sectskipamount=0pt plus30pt
 \def\newsect #1\par{\displayno=0 
   \advance\sectno by 1
   \vskip\sectskipamount\penalty-250\vskip-\sectskipamount
   \bigskip	
   \centerline{\bf \number\sectno. #1}\nobreak
   \medskip	
   \message{#1 }
}
\hyphenation{par-a-met-ri-za-tion par-a-met-ri-zed}
\def\(#1){{\rm(#1)}}\let\leftp=(
\def\activeleftp{\catcode`\(=\active}
{\activeleftp\gdef({\ifmmode\let\next=\leftp \else\let\next=\(\fi\next}}

\def\artkey #1 {{\bf (\number\sectno.#1)}\enspace}

 \let\rmk=\rem

\def\proclaim#1 #2 {\medbreak{\bf#1 }\artkey#2 \bgroup\it\activeleftp}
\def\endproclaim{\egroup\medbreak}
\def\pf{\endproclaim{\bf Proof.}\enspace}
 \let\prp=\prop
\def\lem{\proclaim Lemma }
\def\thm{\proclaim Theorem }

\newcount\displayno
\def\eqlt#1${\global\advance\displayno by 1
 \expandafter\xdef
  \csname \the\sectno#1\endcsname{\number\displayno}
 \eqno\hbox{\rm(\the\sectno-\number\displayno)}$}
\def\tjitem #1 {\global\advance\displayno by 1
 \expandafter\xdef
  \csname \the\sectno#1\endcsname{\number\displayno}
 \item{\the\sectno-\number\displayno} }
\def\eqalignlt#1\cr{\global\advance\displayno by 1
 \expandafter\xdef
  \csname \the\sectno#1\endcsname{\number\displayno}
 \eqno\hbox{\rm(\the\sectno-\number\displayno)}\cr}
\def\disp#1#2{(#1-\csname#1#2\endcsname)}
\def\Cs#1){\unskip~{\rm(\number\sectno.#1)}}
\def\Cd#1){\unskip~{\rm(\number\sectno-\csname\the\sectno#1\endcsname)}}
\def\Co#1-#2){\unskip~{\rm(#1-\csname#1#2\endcsname)}}

\def\item#1 {\par\indent\indent\indent \hangindent3\parindent
 \llap{\rm (#1)\enspace}\ignorespaces}
 \def\part#1 {\par{\rm (#1)\enspace}\ignorespaces}

 \newif\ifproofing \proofingfalse 
 \newcount\refno	 \refno=0
 \def\MakeKey{\advance\refno by 1 \expandafter\xdef
 	\csname\TheKey\endcsname{{%
	\ifproofing\TheKey\else\number\refno\fi}}\NextKey}
 \def\NextKey#1 {\def\TheKey{#1}\ifx\TheKey\NoKey\let\next\relax
  \else\let\next\MakeKey \fi \next}
 \def\NoKey{*!*}
 \def\RefKeys #1\endRefKeys{\expandafter\NextKey #1 *!* }

\def\UThin{\penalty\@M \thinspace\ignorespaces}
\def\relaxnext@{\let\next\relax}
\def\cite#1{\relaxnext@
 \def\nextiii@##1,##2\end@{\unskip\space{\rm[\SetKey{##1},\let~=\UThin##2]}}%
 \in@,{#1}\ifin@\def\next{\nextiii@#1\end@}\else
 \def\next{\unskip\space{\rm[\SetKey{#1}]}}\fi\next}
\newif\ifin@
\def\in@#1#2{\def\in@@##1#1##2##3\in@@
 {\ifx\in@##2\in@false\else\in@true\fi}%
 \in@@#2#1\in@\in@@} \def\SetKey#1{{\bf\csname#1\endcsname}}

\let\texttilde=\~
\def\~{\ifmmode\let\next=\widetilde\else\let\next=\texttilde\fi\next}

\def\p.{\unskip\space p.\UThin}
\def\pp.{\unskip\space pp.\UThin}

 \font\twelvebf=cmbx12
 \font\smc=cmcsc10

\def\eightpoint{\eightpointfonts
 \setbox\strutbox\hbox{\vrule height7\p@ depth2\p@ width\z@}%
 \eightpointparameters\eightpointfamilies
 \normalbaselines\rm
 }
\def\eightpointparameters{%
 \normalbaselineskip9\p@
 \abovedisplayskip9\p@ plus2.4\p@ minus6.2\p@
 \belowdisplayskip9\p@ plus2.4\p@ minus6.2\p@
 \abovedisplayshortskip\z@ plus2.4\p@
 \belowdisplayshortskip5.6\p@ plus2.4\p@ minus3.2\p@
 \skewchar\eighti='177 \skewchar\sixi='177
 \skewchar\eightsy='60 \skewchar\sixsy='60
 \hyphenchar\eighttt=-1
 }
\newfam\smcfam
\def\eightpointfonts{%
 \font\eightrm=cmr8 \font\sixrm=cmr6
 \font\eightbf=cmbx8 \font\sixbf=cmbx6
 \font\eightit=cmti8
\font\eightsmc=cmcsc8
 \font\eighti=cmmi8 \font\sixi=cmmi6
 \font\eightsy=cmsy8 \font\sixsy=cmsy6
 \font\eightsl=cmsl8 \font\eighttt=cmtt8}
\def\eightpointfamilies{%
 \textfont\z@\eightrm \scriptfont\z@\sixrm  \scriptscriptfont\z@\fiverm
 \textfont\@ne\eighti \scriptfont\@ne\sixi  \scriptscriptfont\@ne\fivei
 \textfont\tw@\eightsy \scriptfont\tw@\sixsy \scriptscriptfont\tw@\fivesy
 \textfont\thr@@\tenex \scriptfont\thr@@\tenex\scriptscriptfont\thr@@\tenex
 \textfont\itfam\eightit	\def\it{\fam\itfam\eightit}%
 \textfont\slfam\eightsl	\def\sl{\fam\slfam\eightsl}%
 \textfont\ttfam\eighttt	\def\tt{\fam\ttfam\eighttt}%
 \textfont\smcfam\eightsmc	\def\smc{\fam\smcfam\eightsmc}%
 \textfont\bffam\eightbf \scriptfont\bffam\sixbf
   \scriptscriptfont\bffam\fivebf	\def\bf{\fam\bffam\eightbf}%
 \def\rm{\fam0\eightrm}%
 }

\catcode`\@=13

\ifx\DoublepageOutput\TRUE \doublepageoutput \fi

 \def\SetRef#1 #2\par{%
   \hang\llap{[\csname#1\endcsname]\enspace}%
   \ignorespaces#2\unskip.\endgraf}
 \newbox\keybox \setbox\keybox=\hbox{[18]\enspace}
 \newdimen\keyindent \keyindent=\wd\keybox
 \def\references{\vskip-\smallskipamount
  \bgroup   \eightpoint   \frenchspacing
   \parindent=\keyindent  \parskip=\smallskipamount
   \everypar={\SetRef}}
 \def\endreferences{\egroup}

  \def\paper{\unskip, \bgroup\it}
 \def\paperinfo{\unskip, \egroup}
 \def\preprint{\unskip, \egroup Preprint }
 \def\ag#1 {\preprint alg-geom/#1}
 \def\book{\unskip, ``}

 \def\bookinfo#1{\unskip," #1}

 \def\serial#1#2{\expandafter\def\csname#1\endcsname ##1 ##2 ##3
  {\unskip, \egroup #2 {\bf##1} (##2), ##3}}
 \serial{acta}{Acta Math.}
 \serial{crasp}{C. R. Acad. Sci. Paris}
 \serial{comp}{Comp. Math.}
 \serial{pmihes}{Publ. Math. I.H.E.S.}
 \serial{invent}{Invent. Math.}
 \serial{ma}{Math. Ann.}
 \serial{mathz}{Math. Z.}
 \serial{nagoya}{Nagoya Math. J.} \let\nmj=\nagoya
 \datestrue			

\RefKeys
 BE C83 C86 GLP H Hirsh JK Katz Mori Nj Ramella Verdier
 \endRefKeys

\def\author{T. Johnsen and S. Kleiman}
\def\shorttitle{Toward Clemens' conjecture}

\long\def\topstuff{
\null\vskip24pt plus 12pt minus 12pt
\twelvebf			
\centerline{TOWARD CLEMENS' CONJECTURE}\medskip
\centerline{in degrees between 10 and 24}

 \rm				

 \footnote{}{\noindent 1980 {\it Mathematics Subject Classification}
   (1985  {\it Revision}).  Primary 14J30; Secondary 14H45, 14N10.}

 \vskip12pt plus6pt minus3pt
 \centerline{\smc
 Trygve JOHNSEN\footnote{$^{1}$}{%
    Supported in part by the Norwegian Research Council for Science and
    the Humanities.  It is a pleasure for this author to thank the
Department of Mathematics of the University of Sofia for organizing the
remarkable conference in Zlatograd during the period August 28-September 2,
1995.
It is also a pleasure to thank the
M.I.T. Department of Mathematics for its hospitality from January 1 to
July 31, 1993, when  this work was started.}
 and Steven L.~KLEIMAN,\footnote{$^{2}$}{%
    Supported in part by NSF grant 9400918-DMS.}}

\vskip15pt plus12pt minus12pt
{\parindent=24pt \narrower \noindent \eightpoint
 {\smc Abstract.}\enspace
 We introduce and study a likely condition that implies the following
form of Clemens' conjecture in degrees $d$ between 10 and 24: given a
general quintic threefold $F$ in complex $\IP^4$, the Hilbert scheme of
rational, smooth and irreducible curves $C$ of degree $d$ on $F$ is
finite, nonempty, and reduced; moreover, each $C$ is embedded in $F$
with balanced normal sheaf $\O(-1)\oplus\O(-1)$, and in $\IP^4$ with
maximal rank.
 \par }
 } 

 \topstuff

\def\today{alg-geom/9601024}	

\newsect	Introduction

 Ten years ago, Clemens posed a conjecture about the rational curves on
a general quintic threefold $F$ in complex $\IP^4$.  At once, S.~Katz
\cite{Katz} considered the conjecture in the following form: {\it the
Hilbert scheme of rational, smooth and irreducible curves $C$ of degree
$d$ on $F$ is finite, nonempty and reduced; in fact, each curve is
embedded with balanced normal sheaf $\O(-1)\oplus\O(-1)$}.  Katz
proved this statement for $d\le7$.  Recently, Nijsse \cite{Nj} and the
authors \cite{JK} independently proved the statement for $d\le9$ by
developing Katz's argument.  In the present paper, we will discuss the
possibility of developing Katz's approach further, especially in the
range $10\le d\le 24$.  Notably, we'll focus on a likely condition on a
certain closed subset $I_d'$ of the incidence scheme $I_d$ of all
$C$ and $F$.  In Section~2, we'll derive some consequences from the
condition, including the above form of Clemens' conjecture for
$d\le24$.  In Section~3, we'll discuss some evidence supporting the
condition.

For $d\le9$, a stronger statement holds: the incidence scheme $I_d$ is
reduced and irreducible of dimension 125.  In fact, Katz proved that,
for any $d$, if $I_d$ is irreducible, then the above form of Clemens'
conjecture is true.  Katz established the irreducibility of $I_d$ for
$d\le7$, and Nijsse and the authors established it, via different
arguments, for $d=8,9$.  Moreover, when $I_d$ is irreducible, then, on
a general $F$, each $C$ has several significant additional properties;
see \cite{JK, Cor.~(2.5)}.  For example, each $C$ has maximal rank in
$\IP^4$; that is, for every $k\ge1$, the natural restriction map,
	$$H^0(\IP^4,\O_{\IP^4}(k))\to H^0(C,\O_C(k)),$$
 is either surjective or injective (or both).  These additional
properties play an important role in the authors' work in \cite{JK} on
Clemens' full conjecture, which is discussed briefly below.

On the other hand, $I_d$ is reducible for $d\ge12$, according to
Proposition~(3.2) below.  In fact, $I_d$ always has at least one
component of dimension 125 dominating the space $\IP^{125}$ of all
quintic threefolds $F$; see Lemma~(2.4).  This component was
constructed, more or less explicitly, for infinitely many $d$ by
Clemens \cite{C83, Thm.~(1.27), \p.26}, and for all $d$ by Katz
\cite{Katz, \p.153} (who observed that the general case follows via
Clemens' deformation-theoretic argument from an existence result of
Mori's).  We will see in Section~3 that $I_d$ contains some special
subsets, which do not dominate $\IP^{125}$.  One of them has dimension
$2d+101$ for $d\ge10$, so yields a second component of $I_d$ for
$d\ge12$.  For $d=10,11$ it is unclear whether $I_d$ is irreducible or
not.

Clemens' conjecture is, of course, no less likely to be true.  In fact,
in the above form, it is implied, for $d\le24$, along with the other
consequences of irreducibility, by a likely weaker condition.  This
condition concerns another component of $I_d$, which exists when
$d\le24$.  We'll call it the {\it principal component}, and denote it
by $\?I_{d,0}$.  It is constructed as follows.  In the space $M_d$ of
all $C$, form the open subset $M_{d,0}$ where $H^1(\I_C(5))$ vanishes;
here $\I_C$ denotes the ideal of $C$ in $\IP^4$.  Form the preimage
$I_{d,0}$ in $I_d$ of $M_{d,0}$.  Then $\?I_{d,0}$ is simply the
closure of $I_{d,0}$.  For $d\le24$, we expect that $\?I_{d,0}$ is
equal to the Clemens--Katz component.  In fact, we expect that
$\?I_{d,0}$ is the only component of $I_d$ that dominates $\IP^{125}$.
In other words, we expect that, if $I_d':=I_d-I_{d,0}$, then {\it
$I_d'$ does not dominate $\IP^{125}$}.  This, finally, is our proposed
weaker condition for $d\le24$.  One reasonable way to try to prove it
is to look for a decomposition of $I_d'$ into manageable pieces, each
of which can be shown not to dominate $\IP^{125}$.  On the other hand,
for $d\ge25$, the preimage $I_{d,0}$ is empty, and so the geometry of
$I_d$ is radically different in this range.

Clemens \cite{C86, \p.639} strengthened his conjecture, after Katz's
work, by adding these two assertions: {\it all the rational, reduced
and irreducible curves on a general $F$ are smooth and mutually
disjoint; and the number $n_d$ of curves of degree $d$ is divisible by
$5^3\cdot d$}.  These assertions are not completely true.  Vainsencher
proved that $F$ contains 17,601,000 six-nodal quintic plane curves.
Ellingsrud and Str\o mme and, independently, Candelas, De la Ossa,
Green, and Parkes found that $n_3$ is equal to 371,206,375, which is
divisible by $5^3$, but not 3.  In fact, in their landmark work
introducing mirror symmetry, Candelas et.\ al.\ developed an algorithm
that produces, for any $d$, a number, which they conjecture is equal to
$n_d$.  Afterwards, Kontsevich gave a somewhat different algorithm,
which, he conjectured, also gives the $n_d$.  Although Kontsevich too
was inspired by mathematical physics, his treatment is more
algebraic-geometric.  Moreover, its numbers clearly count both smooth
and nodal curves, which must be connected, but may be reducible.
However, the authors \cite{JK, Thms.~(3.1) and (4.1)} proved that the
only singular, reduced and irreducible, rational curve of degree at
most $9$ on $F$ is a six-nodal plane quintic and that there is on $F$
no pair of intersecting rational, reduced and irreducible curves whose
degrees total at most 9; thus the enumerative significance of
Kontsevich's numbers is established in degree at most $9$.  The
case of degree at least 10 remains open.

Throughout this paper, we use the following notation, which has
already been introduced informally:\smallbreak
 \item a $M_d$ denotes the open subscheme of the Hilbert scheme of
$\IP^4$ parametrizing the rational, smooth and irreducible curves $C$
of degree $d$;
 \item b $\IP^{125}$ denotes the projective space parametrizing the
quintic threefolds $F$ in $\IP^4$;
 \item c $I_d$ denotes the ``incidence'' subscheme of $M_d\times
\IP^{125}$ of pairs $(C,F)$ such that $C\subset F$;
 \item d $M_{d,0}$ denotes the subset of $M_d$ parametrizing the curves
$C$ such that $h^1(\I_C(5))=0$ where $\I_C$ denotes the ideal of $C$ in
$\IP^4$;
 \item e  $I_{d,0}$ denotes the preimage in $I_d$ of $M_{d,0}$;
 \item f  $\?I_{d,0}$ denotes the closure of $I_{d,0}$;
 \item g  $I_d'$ denotes the complement, $I_d-I_{d,0}$.

\newsect	The principal component

In this section, we'll derive some consequences from the (likely)
condition that the (closed) set $I_d'$ does not dominate the space
$\IP^{125}$ of quintic threefolds.  Our main result is Theorem~(2.7);
it asserts that this condition implies Katz's form of Clemens'
conjecture.  The theorem will be derived from Proposition~(2.5), which
asserts this: if $I_d'$ doesn't dominate $\IP^{125}$, then its
complement $I_{d,0}$ does; in fact, then the closure $\?I_{d,0}$ is the
one and only component that does.  We'll call $\?I_{d,0}$ the {\it
principal component\/} of $I_d$.  We'll also prove Proposition~(2.2),
which asserts that, if $I_d'$ doesn't dominate $\IP^{125}$, then, on a
general quintic threefold $F$, the rational curves $C$ possess certain
significant properties; for example, each $C$ has maximal rank in
$\IP^4$.

\lem1 If $d\le24$, then $I_{d,0}$ is smooth, irreducible, and of
dimension $125$; moreover, it dominates $M_{d,0}$, it's open in $I_d$,
and its closure $\?I_{d,0}$ is a component. If $d\ge25$, then $I_{d,0}$
is empty.
 \pf
 It is well known that $M_d$ is smooth of dimension $5d+1$ and is
irreducible.  Moreover, these properties are not hard to establish.
Indeed, fix $C\in M_{d}$.  The restricted Euler sequence,
	$$0\to\O_C\to\O_C(1)^{\oplus5}\to{\cal T}_{\IP^4}|C\to0,\eqlt2$$
 yields $H^1({\cal T}_{\IP^4}|C)=0$.  So the sequence of tangent and
normal sheaves,
 $$0\to{\cal T}_C\to{\cal T}_{\IP^4}|C\to{\cal N}_{C/\IP^4}\to0,\eqlt3$$
 yields $H^1({\cal N}_{C/\IP^4})=0$.  Hence, by the standard theory of
the Hilbert scheme, $M_d$ is smooth at $C$ of dimension $h^0({\cal
N}_{C/\IP^4})$, and the latter number can be found  using the same
two exact sequences.  Finally, $M_d$ is irreducible as it's the image
of an open subset of the space of parametrized maps from $\IP^1$ to
$\IP^4$, and this space of maps is just the space of 5-tuples of
homogeneous polynomials of degree $d$ in two variables.

Again, fix $C\in M_{d}$.  Then $C\in M_{d,0}$ if and only if the
natural map,
	$$H^0(\IP^4,\O_{\IP^4}(5))\to H^0(C,\O_C(5)),\eqlt1$$
 is surjective, thanks to the long exact cohomology sequence.  Hence,
if $d\le7$, then $M_{d,0}=M_d$.  Indeed, the surjectivity of \Cd1) is
obvious if $C$ is a line, a conic, or a twisted cubic.  If $4\le
d\le7$, then $C$ cannot lie in plane, and the surjectivity holds by the
theorem on \p.492 of \cite{GLP}.  If $8\le d\le25$, then $M_{d,0}$ is
nonempty; indeed, if $C\in M_d$ is general, then the surjectivity of
\Cd1) holds by the maximal-rank theorem \cite{BE, Thm.~1, \p.215},
because the source and target have dimensions 126 and $5d+1$.  If
$d\ge26$, then surjectivity cannot hold, and so $M_{d,0}$ is empty.

The subset $M_{d,0}$ of $M_d$ is open for any $d$.  Indeed, let ${\bf
C}$ be the universal curve in $\IP^4\x M_d$, and $\I_{\bf C}$ its ideal.
Then $\I_{\bf C}$ is flat over $M_d$.  Hence $h^1(\I_{\bf C}(5))$  is
upper semi-continuous.  Therefore, $M_{d,0}$ is open.  Hence, its
preimage $I_{d,0}$ is open in $I_d$, and its closure will be a
component provided  $I_{d,0}$ is nonempty and irreducible.

Let $C\in M_{d,0}$.  Then, by definition, $H^1(\I_C(5))$ vanishes.
Hence the direct image ${\cal Q}$ of $\I_{\bf C}(5)$ is locally free on
$M_{d,0}$, and its formation commutes with base change to the fibers.
Hence $I_{d,0}$ is equal to $\IP({\cal Q}^*|M_{d,0})$.  Now, for
$d\le25$, since \Cd1) is surjective, $H^0(\I_C(5))$ has dimension
$125-5d$.  Hence, $H^0(\I_C(5))$ is zero for $d=25$, and is nonzero for
$d\le24$.  For $d\le24$, therefore, $I_{d,0}$ is smooth, irreducible,
of dimension $125$, and dominates $M_{d,0}$.  Thus the lemma is proved.

\prp2 Assume that $I_d'$ does not dominate $\IP^{125}$.  Let $F$ be a
general quintic threefold in $\IP^4$, and let $C$ be a rational, smooth
and irreducible curve of degree $d$ at most $24$ on $F$.
 \part 1 Then $C$ is embedded in $\IP^4$ with maximal rank.
 \part 2 Form the restriction to $C$ of the twisted sheaf of
differentials of $\IP^4$.  Then this locally free sheaf has\/ {\rm generic}
splitting type; namely, if $d=4n+m$ where $0\le m<4$, then
	$$\Omega^1_{\IP^4}(1)|C =\O_C(-n-1)^m\oplus\O_C(-n)^{4-m}.$$
 \part 3 If $d\le4$, then $C$ is a rational normal curve of degree $d$,
and if $d\ge4$, then $C$ spans $\IP^4$.
 \part 4 If $d=1$, then $C$ is $1$-regular; if $2\le d\le4$, then $C$
is $2$-regular; if $5\le d\le7$, then $C$ is $3$-regular; if $8\le
d\le11$, then $C$ is $4$-regular; if $15\le d\le17$, then $C$ is
$5$-regular; and if $18\le d\le24$, then $C$ is $6$-regular.
 \pf This result was proved in \cite{JK, Cor.~(2.5)} (without the
hypothesis on $I_d'$) for $d\le9$.  For $10\le d\le24$, the proof is
similar.  First, observe that, since $F$ is general, $C$ does not lie
in any given proper closed subset $N$  of $M_d$.  Indeed, the preimage of
$N$ in $I_d$ consists of two parts, the part in $I_d'$ and that in
$I_{d,0}$.   Neither part dominates $\IP^{125}$: the first doesn't
by hypothesis, and the second doesn't by virtue of Lemma~\Cs1), which
implies that this part has dimension at most $124$.

To prove (2), apply the observation above to the subset $N$ of $M_d$ of
curves without the asserted splitting type; $N$ is a proper closed
subset by a theorem of Verdier's \cite{Verdier, Thm., p.139} (see also
\cite{Ramella, Thm.~1, p.~181}).  To prove (3), apply the observation
above to the subset $N$ of $M_d$ of curves not spanning $\IP^4$; here
$N$ is proper if $d\ge4$, because, clearly, $\dim N\le4d+4$ whereas
$\dim M_d=5d+1$.  To prove (1), apply the observation above to the
subset $N$ of $M_d$ of curves that either don't span $\IP^4$ or aren't
of maximal rank; here $N$ is proper if $d\ge4$ by (3) and by the
maximal-rank theorem \cite{BE, Thm.~1, \p.215}.  Finally, (1) implies
(4) by virtue of the long exact sequence of cohomology extending the
map~\Cd1).

\lem3 Let $(C,F)\in I_d$, and assume that $F$ is smooth along $C$.
Then the following conditions are equivalent:\smallskip
 \item i At $(C,F)$, the incidence scheme $I_d$ is smooth of dimension
$125$, and the differential $d\beta$ of the projection
$\beta\:I_d\to\IP^{125}$ is surjective.
 \item ii At $C$, the Hilbert scheme of $F$ is reduced of dimension $0$.
 \item iii The normal sheaf $\N_{C/F}$ has a balanced decomposition,
	     $$\N_{C/F}=\O_{\IP^1}(-1)\oplus\O_{\IP^1}(-1).$$
 \smallskip\noindent
 If any one of these conditions obtains, then $(C,F)$ lies on a unique
component of $I_d$, which is generically reduced, has dimension $125$,
and dominates $\IP^{125}$.
 \pf It is necessary and sufficient for (i) to hold that, at $(C,F)$,
the fiber of $\beta$ be smooth of dimension 0 and that $\beta$ be flat.
However, in any event, $I_d$ is simply an open subscheme of the
relative Hilbert scheme $\Hilb_{{\bf F}/\IP^{125}}$ where ${\bf F}$ is
the universal family of quintics.  Hence (i) implies (ii).  Moreover,
(i) is implied by (iii), because, by the standard theory of the
relative Hilbert scheme, when $H^1(\N_{C/F})$ vanishes, then
$\Hilb_{{\bf F}/\IP^{125}}$ is smooth over $\IP^{125}$ with
$H^0({\cal N}_{C/F})$ as fiber dimension.

It is also part of standard theory that $H^0(\N_{C/F})$ is equal to the
Zariski tangent space to the Hilbert scheme of $F$ at the point $C$;
hence (ii) holds if and only if $H^0(\N_{C/F})$ vanishes.  Now, it is
easy to see that $\N_{C/F}$ has as determinant $\O_{\IP^1}(-2)$.
Indeed,  the sequences \Cd2) and \Cd3) show
that $\N_{C/\IP^4}$ has as determinant $\O_{\IP^1}(5d-2)$.  Then the
sequence of normal sheaves,
	$$0\to\N_{F/\IP^4}\to\N_{C/\IP^4}\to\N_{C/F}\to0$$
 yields the determinant of $\N_{C/F}$.  Now,
$\N_{C/F}=\O_{\IP^1}(a)\oplus\O_{\IP^1}(b)$ for some $a$ and $b$. Hence
$a+b=-2$.  Hence $H^0(\N_{C/F})$ vanishes if and only if both $a$ and
$b$ are $-1$.  Therefore, (ii) and (iii) are equivalent to each other,
whence also to (i).

Suppose one of the conditions (i), (ii) or (iii) obtains; then all
three do.  Hence (i) implies that $(C,F)$ lies in the smooth locus of
$I_d$, so in a unique component, which is reduced at $(C,F)$.
Moreover, (i) implies that this component has dimension $125$, and that
the projection $\beta$ onto $\IP^{125}$ is smooth at $(C,F)$, so open
on a neighborhood of it.  Therefore, the component of $(C,F)$ dominates
$\IP^{125}$, and the proof is complete.

\lem4 The incidence scheme $I_d$ always has at least one component that
is generically reduced, that has dimension $125$, and that dominates
$\IP^{125}$.
 \pf By the work of Clemens and Katz (see \cite{Katz, Thm.~2.1,
p.~153}), there is a pair $(C,F)\in I_d$ such that $F$ is smooth along
$C$ and such that the normal sheaf $\N_{C/F}$ has a balanced
decomposition, $\O(-1)\oplus\O(-1)$.  Hence, Lemma~\Cs3) yields the
assertion.

\prp5 Assume that $I_d'$ does not dominate $\IP^{125}$.  Then $d\le24$,
and the principal component $\?I_{d,0}$ is the one and only component
of $I_d$ that dominates $\IP^{125}$.
 \pf By Lemma~\Cs4), there is at least one component of $I_d$ that
dominates $\IP^{125}$.  Given any such component, it cannot lie in
$I_d'$ by hypothesis; so it lies in the closure of the complement of
$I_d'$, namely, $\?I_{d,0}$.  So $\?I_{d,0}$ is nonempty.  Hence
Lemma~\Cs1) implies that $d\le24$ and that $\?I_{d,0}$ is a component.
The remaining assertions now follow.

\lem6 Let $\~I_d$ be a component of $I_d$, and assume that $\~I_d$ is
generically reduced, has dimension $125$, and dominates $\IP^{125}$.
Let $F\in\IP^{125}$ be a general quintic, and let $\Phi$ be the set of
$C$ with $(C,F) \in \~I_d$.  Then $\Phi$ is finite and nonempty.
Moreover, at each $C$ in $\Phi$, the Hilbert scheme of $F$ is reduced
of dimension $0$; in fact, each $C$ is embedded in $F$ with balanced
normal sheaf, $\O(-1)\oplus\O(-1)$.
 \pf The set $\{(C,F)|C\in\Phi\}$ is simply the fiber of $\~I_d$ over
$F$.  So it is finite and nonempty, because $\~I_d$ has dimension $125$
and dominates $\IP^{125}$ and because $F$ is general.  By the same
token, this fiber lies in the smooth locus of $\~I_d$, which is
nonempty because $\~I_d$ is generically reduced.  Hence, by Sard's
lemma, the differential of the projection $I_d\to\IP^{125}$ is
surjective along $\Phi$.  Therefore, the remaining assertions follow
from Lemma~\Cs3).

\thm7 Assume that $I_d'$ does not dominate $\IP^{125}$.  Then $d\le24$.
Let $F$ be a general quintic threefold in $\IP^4$, and in the Hilbert
scheme of $F$, form the open subscheme of rational, smooth and
irreducible curves $C$ of degree $d$.  Then this subscheme is finite,
nonempty, and reduced; in fact, each $C$ is embedded in $F$ with
balanced normal sheaf $\O_{\IP^1}(-1)\oplus\O_{\IP^1}(-1)$, and
possesses the properties (1) to (4) of Proposition~(2.2).
 \pf Proposition~\Cs2) applies, so its properties (1) to (4) hold (but,
so far, possibly are vacuous).  By Proposition~\Cs5) above, $d\le24$,
and $\?I_{d,0}$ is the one and only component of $I_d$ that dominates
$\IP^{125}$.  This component is generically reduced and has dimension
125 by Lemma~\Cs1), or alternatively by Lemma~\Cs4).  Hence Lemma~\Cs6)
yields the remaining assertions.

\newsect	Other subsets

In this section, we'll prove Proposition \Cs2), which asserts that
$I_d$ is reducible for $d\ge12$.  We'll proceed by introducing and
studying some basic subsets $J^e_d$ and $K_d$ of $I_d$.  For $d\ge12$,
they provide one or more components of $I_d$, which do not dominate
$\IP^{125}$.  After proving the proposition, we'll make two remarks;
the first discusses a refinement of the condition that $I_d'$ does not
dominate  $\IP^{125}$, and the second discusses the location in $I_d$ of
the pair $(C,F)$ of Clemens and Katz.

The subsets $J^e_d$ and $K_d$ of $I_d$ are the following: \smallbreak
 \item a $J_d^e$ is the set of pairs $(C,F)\in I_d$ such that
$C$ spans a hyperplane $H$ and lies on a smooth surface $S$ of
degree $e$ in $H$;
 \item b $K_d$ is the set of pairs $(C,F)\in I_d$ such that $C$ spans a
hyperplane $H$ and $H^0(\I_{C/H}(5))=0$ where $\I_{C/H}$ is the ideal
of $C$ in $H$.

\lem1 The dimensions of the above sets are as follows:
	$$\eqalign{\dim J_d^2=2d+101\for d\ge10;\quad
	&\dim J_d^3=d+101\for d\ge15;\cr
	\dim J_d^4\le97\for d\ge20;\quad
	&\dim K_d=4d+73\for d\ge11.\cr}$$
 None of these sets dominates $\IP^{125}$.  Moreover, $K_d$ is empty
for $d\le10$.
 \pf
 Fix $(C,F)\in J_d^e$.  By definition, $C$ spans a hyperplane $H$ and
lies on some smooth surface $S$ of degree $e$ in $H$.  If $d\ge e^2$,
then $S$ is uniquely determined; otherwise, $C$ would lie in the
intersection of two different smooth surfaces of degree $e$ in $H$, so
$C$ would be equal to this intersection, and so would have nonzero
genus.  Furthermore, if $d\ge5e$, then $S$ lies in $F$; otherwise, the
intersection of $S$ and $F$ would be a curve containing $C$, so $C$
would be equal to this intersection, and so would have nonzero genus.

 Vary $(C,F)\in J_d^e$, and form the space $\~J_d^e$ of corresponding
triples $(C,S,F)$.  If $e\le5$ and $d\ge5e$, then, by the preceding
argument, the projection $\~J_d^e\to J_d^e$ is bijective, so $J_d^e$
and $\~J_d^e$ have the same dimension and the same image in
$\IP^{125}$.  We'll now compute this dimension and image for $e=2,3,4$.
The fiber of $\~J_d^e$ over a pair $(S,F)$ consists of all $C$ in $S$.
So this fiber has dimension $2d-1$ if $e=2$, dimension $d-1$ if $e=3$,
and dimension at most 0 if $e=4$.  These dimensions are well known, and
they are easy to check.  (For $e=3$, use \cite{H, 4.8, p.~401} and
\cite{H, 4.8, p.~407}.  For e=4, note that there are at most finitely
many $C$ on a given $S$ because the normal sheaf of $C$ is equal to
$\O_{\IP^1}(-2)$; in fact, a general $S$ can contain no $C$ because all
the curves on it are complete intersections by Noether's theorem).

The $F$ containing a fixed $S$ form a space of dimension
$h^0(\I_S(5))-1$.  To compute it, use the natural exact sequence of
ideals,
	$$0\to\I_H\to\I_S\to\I_{S/H}\to0.$$
 The first term is equal to $\O_{\IP^4}(-1)$ and the third to
$\O_{\IP^3}(-e)$.  Hence
	$$h^0(\I_S(5))=h^0(\O_{\IP^4}(4))+h^0(\O_{\IP^3}(5-e))
	=70+{8-e\choose3}.$$
 The various $S$ in a fixed  $H$ form a space of
dimension ${3+e\choose3}-1$, and the various $H$ form a $\IP^4$.  Hence
the various pairs $(S,F)$ form a space of dimension,
	$$\eqalign{(70+20-1)+(10-1+4)&=102\hbox{ if } e=2,\cr
	(70+10-1)+(20-1+4)&=102\hbox{ if } e=3,\cr
	(70+4-1)+(35-1+4)&=111\hbox{ if } e=4.\cr}$$
  These numbers are less than 125.  Therefore, $J_d^e$ doesn't dominate
$\IP^{125}$ for $e=2,3,4$ and $d\ge5e$.  Furthermore,
 $$\eqalign{\dim J_d^2&=(2d-1)+102=2d+101\for d\ge10,\cr
	\dim J_d^3&=(d-1)+102=d+101\for d\ge15,\cr
	\dim J_d^4&\le0+111=111\for d\ge20.\cr}$$
 Thus the assertions about the $J_d^e$ are proved.

To analyze $K_d$, observe that, in $M_d$, the $C$ that lie in a fixed
hyperplane $H$ form a closed subset of dimension $4d$, and that, in this
closed subset, those $C$ with $h^0(\I_{C/H}(5))=0$ form an open subset
by upper semi-continuity of dimension.  This open set is nonempty, so of
dimension $4d$, if and only if $d\ge11$; indeed, the maximal rank
theorem for rational curves in $\IP^3$ \cite{Hirsh, Thm.~0.1, p.~209}
implies that, for a general $C$ in $H$, the natural map,
	$$H^0(H,\O_H(5))\to H^0(C,\O_C(5)),$$
 is injective if and only if $d\ge11$, because the source and target
have dimensions 56 and $5d+1$.  Hence $K_d$ is empty for $d\le10$.  On
the other hand, since the various $H$ form a $\IP^4$, the image of
$K_d$ in $M_d$ therefore has dimension $4d+4$ for $d\ge11$.

Whenever $C\subset H$ and $H^0(\I_{C/H}(5))=0$, the natural inclusion map,
	$$H^0(\I_H(5)) \to H^0(\I_C(5)),\eqlt2$$
 is bijective.  Since the source has dimension 70, the fiber in $K_d$
over $C$ is a $\IP^{69}$.  Hence $K_d$ has dimension $(4d+4)+69$, or
$4d+73$.  Moreover, since \Cd2) is bijective, the image of $K_d$ in
$\IP^{125}$ is equal to the set of quintics $F$ that contain a
hyperplane.  The latter set has dimension $69+4$, or $73$.  The proof
is now complete.

\prp2 If $d\ge12$, then $I_d$ is reducible.  In fact, if $d\ge13$, then
$I_d$ has a component of dimension at least $126$, as well as one of
dimension  $125$.
 \pf On the one hand, $I_d$ always has at least one component that has
dimension $125$ and that dominates $\IP^{125}$ by Lemma~(2.4).  On the
other hand, if $d\ge10$, then $I_d$ has a subset, namely $J_d^2$, that
has dimension $2d+101$ and that doesn't dominate $\IP^{125}$ by
Lemma~(3.1).  Hence, if $d\ge13$, then $\dim J_d^2\ge 126$, and so $I_d$
has a component of dimension at least 126, as well as one of dimension
$125$.  Suppose $d=12$.  Then $J_d^2$ has dimension 125, but doesn't
dominate $\IP^{125}$.  So $J_d^2$  cannot lie in the component of $I_d$
that dominates $\IP^{125}$.  Hence $I_d$ is still reducible.  Thus the
proposition is proved.

\rmk3  For $d\le24$, it is not unreasonable to hope that the complement
$I_d''$ of $\?I_{d,0}$ in $I_d$ lies in the closure of the union of
$J_d^2$, $J_d^3$, and $K_d$, and that this union doesn't dominate
$\IP^{125}$.  If this hope is confirmed, then $I_d'$ doesn't dominate
$\IP^{125}$ either, because $I_d'-I_d''$ is equal to $\?I_{d,0} -
I_{d,0}$ and so has dimension at most 124.  Hence, then the conclusions
of Proposition~(2.2), Proposition~(2.5), and Theorem~(2.7) will hold.

Lemma~\Cs1) supports this hope.  Indeed, the lemma implies that $K_d$
for $d\ge13$ yields another component of $I_d$ that doesn't dominate
$\IP^{125}$, and that $J_d^3$ for $d\ge24$ yields one too, but that
$J_d^4$ for $d\ge20$ does not.  Of course, to confirm our hope, we must
handle the $J_d^e$ for the $d$ and $e$ not covered by Lemma~\Cs1), and
we must handle the subset of $I_d$ of pairs $(C,F)$ such that $C$ spans
a hyperplane $H$ and lies on a singular, reduced and irreducible,
surface of degree $e$ in $H$, but not on a smooth one.  However, we may
assume that $e\le5$, because $C$ lies in the intersection of $H$ and
$F$, and the latter will be a surface of degree 5 for a suitable $F$ if
$(C,F)\notin K_d$.  Of course, we may assume $d\ge10$ because $I_d$ is
irreducible for $d\le9$ by \cite{JK, Prp.~(2.2)}.  Moreover, we may
assume that $e\ge3$ because, if $C$ lies in a plane, then $d$ is 1 or
2, and if $C$ lies on a 2-dimensional singular quadric cone, then
$d\le3$ by \cite{H, Ex.~2.9, p.~384}.

To confirm our hope, we must also handle the subset of $I_d$ of pairs
$(C,F)$ such that $C$ spans $\IP^4$ and lies on a hypersurface $T$ of
degree $t$ with $2\le t\le5$.  Now, for $t=2,3,4$, this subset does not
trivially yield a new component of $I_d$.  Indeed, locally $dt+1$
conditions must be satisfied for a $C$ in $M_d$ to lie on a given $T$,
and each such $C$ lies at least in the reducible quintics $F$ of the
form $T+T'$.  Hence, the various such $(C,F)$ form a space of dimension
at least,
	$$(5d+1)-(dt+1)+{4+t\choose4}-1+{9-t\choose4}-1.$$
 This number is equal to $3d+48$ for $t=2$, to $2d+48$ for $t=3$, and to
$d+73$ for $t=4$. Its maximum value is achieved for $d=24$ and $t=2$,
and this maximum is only 120, not enough to yield a new component.

It is less likely (as well as unnecessary) that $I_d'$ lies in the
closure of the union of $J_d^2$, $J_d^3$, and $K_d$.  In other words,
there may be pairs $(C,F)$ outside this closure, yet in $\?I_{d,0} -
I_{d,0}$.  For example, such a pair might arise from a curve $C$ of
degree $9$ that spans $\IP^4$ and has a 7-secant.

\rmk4 It is interesting to look at the pair $(C,F)$ found by Katz
\cite{Katz, p.~153}, and observe where it sits in $I_d$.  Katz began
with the curve $C\in M_d$ constructed by Mori \cite{Mori}.  It lies on
a smooth quartic surface $S$ in a hyperplane $H$ in $\IP^4$.  So it
lies in all the reducible quintic surfaces $S+L$ where $L$ is a plane
in $H$.  Hence $h^0(\I_{C/H}(5))\ge4$.  So $(C,F)\notin K_d$.
Moreover, if $d\ge12$, then $C$ cannot lie on a cubic surface
(otherwise it would lie on the intersection of this cubic with $S$),
and so $(C,F)\notin J_d^3$.  Similarly, if $d\ge8$, then $(C,F)\notin
J_d^2$.

If $d\ge10$, then $(C,F)\notin I_{d,0}$.  Indeed, $H^1(\I_C(5))=
H^1(\I_{C/H}(5))$ because of the exact sequence of twisted ideals,
	 $$0\to\I_H(5)\to\I_C(5)\to\I_{C/H}(5)\to0,$$
 whose first term is equal to $\O_{\IP^4}(4)$.  Hence it's enough to
check that $h^1(\I_{C/H}(5))>0$.  Now, the usual long exact cohomology
sequence yields
 $$\eqalign{h^1(\I_{C/H}(5))&=h^0(\I_{C/H}(5))-h^0(\O_H(5))+h^0(\O_C(5))\cr
			&\ge(5d+1)-56+4=5d-51.\cr}$$
 Hence $h^1(\I_{C/H}(5))>0$ if $d\ge11$.  A more sophisticated, but
well-known, argument works for $d\ge10$.  Namely, $C$ has a
$(d-3)$-secant line; it's the curve $D$ in \cite{Mori, p.~129}.  By
Bezout's theorem, $D$ lies in every hypersurface of degree $d-4$
containing $C$.  So $C$ is not cut out by such hypersurfaces.  Hence,
$C$ is $(d-4)$-irregular.  Therefore, $H^1(\I_C(d-5))$ is nonvanishing
since $H^q(\I_C(d-4-q))$ vanishes for $q\ge2$.  It follows that
$H^1(\I_C(5))$ is nonvanishing if $d\ge10$.  Thus there's some content
to our conjecture that $(C,F)\in\?I_{d,0}$ for $d\le24$.

\newsect References

\references

BE
E. Ballico and P. Ellia
 \paper On the Postulation of Curves in $\IP^4$
 \mathz 188 1985 215--23

C83
H. Clemens
 \paper Homological equivalence, modulo algebraic equivalence, is not
fin\-ite\-ly generated
 \pmihes 58 1983 19--38

C86
H. Clemens
 \paper Curves on higher-dimensional complex projective manifolds
  \paperinfo Proc. International Cong. Math., Berkeley, 1986, 634--40

GLP
L. Gruson, R. Lazarsfeld, and C. Peskine
 \paper On a theorem of Castelnuovo and the equations defining space
curves
 \invent 72 1983 491--506

H
R. Hartshorne
 \book Algebraic Geometry
 \bookinfo GTM {\bf 52}, Springer-Verlag, 1977

Hirsh
A. Hirschowitz
 \paper Sur la postulation generique des courbes rationelles
 \acta 146 1981 209--30

JK
T. Johnsen and S. Kleiman
 \paper Rational curves of degree at most 9 on a general quintic
threefold \preprint alg-geom/9510015

Katz
S. Katz
 \paper On the finiteness of rational curves on quintic threefolds
 \comp 60 1986 151--62

Mori
S. Mori
 \paper On degrees and genera of curves on smooth quartic surfaces in
$\IP^3$ \nmj 96 1984 127--32

Nj
P. Nijsse
 \paper Clemens' conjecture for octic and nonic curves
 \paperinfo Preprint, U. Leiden, 1993, to appear in Indag. Math

Ramella
L. Ramella
 \paper La stratification du sch\'ema de Hilbert des courbes rationelles
de $\IP^n$ par le fibr\'e tangent restreint
 \crasp 311 1990 181--84

Verdier J.-L. Verdier
 \paper Two dimensional sigma-models and harmonic maps from ${\bf S}^2$
to ${\bf S}^{2n}$
 \paperinfo in ``Group Theoretical Methods in Physics $\bullet$
Proceedings, Istanbul, Turkey,'' Springer Lecture Notes in Physics {\bf
180} (1983), 136--41

\endreferences

  \bigskip

 \eightpoint\smc Mathematical Institute, University of Bergen, All\'egaten
55, N-5007 Berg\-en, Norway
 \medskip
 Department of Mathematics, 2--278 MIT, Cambridge, MA 02139, USA

\bye